# On the structure of completely inverse AG**-groupoids

R. A. R. Monzo[1]

**Abstract.** The main result of this paper is the determination of the structure of completely inverse AG**-groupoids modulo semilattices of abelian groups and their involutive, idempotent-fixed automorphisms. In addition we prove that the class of completely inverse AG**-groupoids, the class of strongly regular AG**-groupoids, and the class of strongly regular AG-groupoids with commuting idempotents are identical. Also, we prove that a groupoid is an AG-group if and only if it is a left simple, strongly regular AG**-groupoid if and only if it is a left simple AG-groupoid with a left identity. Finally, it is proved that a medial groupoid $S$ is an inflation of a completely inverse AG**-groupoid if and only if $S^2$ is a completely inverse AG**-groupoid.

## 1 Introduction

Completely inverse AG**-groupoids are groupoids that satisfy the equations $(xy)z = (zy)x$ and $x(yz) = y(xz)$ and in which each element has a unique inverse with which it commutes. Such groupoids have been shown to have an intricate and plentiful structure [3,4]. This paper explores the intimate connection between completely inverse AG**-groupoids and semilattices of abelian groups.

Our main result is Corollary 22, which states that a groupoid $S$ is a completely inverse AG**-groupoid if and only if there is a semilattice of abelian groups $(S, •)$ and an involutive, idempotent-fixed automorphism $A$ on $S$ such that $ab = Aa • b$, for $a,b \in S$. We also prove in Theorem 3 that the class of completely inverse AG**-groupoids coincides with the class of strongly regular AG**-groupoids and the class of strongly regular AG-groupoids with commuting idempotents. The class of AG-groups is proved in Theorem 8 and Corollary 9 to be equal to the class of left-simple, completely inverse AG**-groupoids and the class of left-simple, AG-groupoids with a left identity.

In semigroup theory it is well-known that a semigroup $S$ is an inflation of a semilattice of groups if and only if $S^2$ is a semilattice of groups. We prove that if a groupoid $S$ satisfies the medial law then $S$ is an inflation of a completely inverse AG**-groupoid if and only if $S^2$ is a completely inverse AG**-groupoid (Theorem 10).

## 2. Preliminaries

An *AG-groupoid* is a groupoid satisfying the *invertive law* $(xy)z = (zy)x$. Such a groupoid satisfies the *medial law* $(xy)(zw) = (xz)(yw)$. An *AG**-groupoid* is an AG-groupoid satisfying the identity $x(yz) = y(xz)$ and such a groupoid satisfies the *paramedial law* $(xy)(zw) = (wy)(zx)$.. An AG-groupoid with a left identity element also satisfies the paramedial law and the identity $x(yz) = y(xz)$ and, hence, it is an AG**-groupoid.

In a groupoid $S$ with $a,b \in S$, $b$ is called a *left inverse* of $a$ if $a=(ab)a$. The element $b$ is called an *inverse of a* if $a$ and $b$ are left inverses of each other. A groupoid $S$ is called *an inverse groupoid* [*a completely inverse groupoid*] if every element $a \in S$ has a unique inverse [a unique inverse with which it commutes]. The *set of inverses* of the element $a$ is denoted by $V(a)$. In an AG**-groupoid $V(a) \neq \emptyset$ if and only if $V(a)$ has exactly a single element; that is, inverses, if they exist, are unique [1, remark 3]. In an inverse groupoid $S$ we can denote the unique inverse of $a$ as $a^{-1}$ and so $aa^{-1} = a^{-1}a$ in a completely inverse groupoid.

If $I: S \rightarrow T$ is an isomorphism then we write $S \cong T$ or $I: S \cong T$. If in addition, $S = T$ then we call $I$ an *automorphism*. We call $I$ an *idempotent-fixed* (or *E-fixed*) automorphism if it is an automorphism and is the identity mapping when restricted to $E = E(S) = \{x \in S : xx=x\}$. We call $I$ an *involutive automorphism* if it is an automorphism and $I^2$ is the identity mapping on $S$. We denote the collection of involutive automorphisms on a groupoid $S$ by $AUT^2(S)$ and the collection of involutive, $E(S)$-fixed automorphisms by $AUT^2_e(S)$.

**Definition [5]** An element $a$ of an AG-groupoid $S$ is called a *strongly regular element of $S$* if there exists $x \in S$ such that $a=(ax)a$ and $ax=xa$. An groupoid $S$ is called a *strongly regular groupoid* if all elements of $S$ are strongly regular.

Note that a strongly regular AG-groupoid is therefore an AG-groupoid in which every element has a left inverse with which it commutes. If $x$ is a left inverse of $a$ and $ax=xa$ then
$ax = [(ax)a]x = (xa)(ax) = (ax)(ax) \in E(S)$ and if $S$ is a completely inverse AG**-groupoid then
$aa^{-1} = a^{-1}a \in E(S)$ and $(ab)^{-1} = a^{-1}b^{-1}$. These facts will be used throughout the paper without mention.

**Notation.** $A <_r S$, $A <_l S$ and $A < S$ will denote that $A$ is a *right, left or two-sided ideal* of a groupoid $S$ (that is, $AS \subseteq A$, $SA \subseteq A$ or $AS \cup SA \subseteq A$).

## 3. Strongly regular AG-groupoids

**Lemma 1.** *If $S$ is a strongly regular AG-groupoid then every element of $S$ has an inverse with which it commutes.*

*Proof* Let $a \in S$ and let $a=(ax)a$, with $ax=xa$. Let $y = (xa)x$. Then,
$(ay)a = \{[(ax)a][(xa)x]\}a = \{[(ax)(xa)](ax)\}a = (ax)a = a$. Also,
$ya = [(xa)x]a = (ax)(xa) = xa = ax = [(ax)a][(xa)x] = ay$. Therefore,
$(ya)y = (xa)y = (ya)x = (xa)x = y$. So $a$ and $y$ commute and are left inverses of each other. ∎

**Lemma 2.** *A strongly regular AG-groupoid $S$ is an AG**-groupoid if and only if $E(S)$ is a semilattice.*

*Proof* $(\Rightarrow)$ As pointed out in the introduction, if an AG-groupoid satisfies the identity
$x(yz) = y(xz)$ then it satisfies the identity $(xy)(zw) = (wy)(zx)$. Let $\{e,f\} \subseteq E(S)$.
Then $ef = (ee)(ff) = (fe)(fe) = fe$. Also, $(ef)g = (gf)e = e(gf) = e(fg)$
and so $E(S)$ is a semilattice in an AG**-groupoid.

$(\Leftarrow)$ Assume that $S$ is a strongly regular AG-groupoid and that $E(S)$ is a semilattice. Let
$a,b,c \in S$. Then, by Lemma 1, there exists $d,e,f \in S$ such that $(ad)a = a$, $(be)b = b$,
$(cf)c = c$, $ad, be, cf \in E(S)$ and $a(bc) = [(ad) \bullet (be) \bullet (cf)][a(bc)]$. Using the invertive and medial laws repeatedly and the fact that $E(S)$ is a semilattice,
$$a(bc) = (ad \bullet a)(bc) = (ad \bullet b)(ac) = [(ad \bullet a)(cf \bullet c)] =$$
$$= (ad \bullet b)[(ad \bullet cf)(ac)] = [(ad)(ad \bullet cf)][b(ac)] = (ad \bullet cf)[b(ac)] =$$
$$= [(ad) \bullet (be) \bullet (cf)]\{(ad \bullet cf)[b(ac)]\} = [(be) \bullet (ad) \bullet (cf)]^2 \{(ad \bullet cf)[b(ac)]\} =$$
$$= \{[(be) \bullet (ad) \bullet (cf)](ad \bullet cf)\}\{[(be) \bullet (ad) \bullet (cf)][b(ac)]\} =$$
$$= [(be) \bullet (ad) \bullet (cf)][b(ac)] = b(ac).$$
Therefore, $S$ is an AG**-groupoid. ∎

**Theorem 3.** *The following statements are equivalent in a groupoid S.*

*(1) S is a completely inverse AG\*\*-groupoid.*
*(2) S is a strongly regular AG-groupoid and E = E(S) is a semilattice.*
*(3) S is a strongly regular AG\*\*-groupoid.*

*Proof* (1$\Rightarrow$2) By definition, a completely inverse AG\*\*-groupoid is clearly a strongly regular AG-groupoid. By Lemma 2, *E* is a semilattice.
(2$\Rightarrow$1 and 3) Suppose that *S* is a strongly regular AG-groupoid and *E* is a semilattice. By Lemma 1, every element has an inverse with which it commutes.
By Lemma 2, *S* is an AG\*\*-groupoid. Since in an AG\*\*-groupoid inverses are unique, *S* is a completely inverse AG\*\*-groupoid.
(3$\Rightarrow$2) By Lemma 2, *E* is a semilattice. ∎

## 4. Left simple, strongly regular AG\*\*-groupoids

**Definition [3].** An AG groupoid with left identity *e* in which for every element *a* there exists an element $a^*$ such that $a^* a = e$ is called an AG -group.
In [3] it was proved that $a^*$ is unique and that $aa^* = e$. We will denote $a^*$ by $a^{-1}$.
Note that $a \in V(a^{-1})$ and $a^{-1} \in V(a)$.

**Lemma 4 [3, Proposition 2.2]** *The following statements are equivalent in an AG -groupoid S with left identity element e.*
*(1) S is an AG -group;*
*(2) for every element $a \in S$ there exists an element $a^*$ such that $aa^* = e$ ;*
*(3) every element $a \in S$ has a unique inverse $a^{-1}$ and*
*(4) the equation xa = b has a unique solution for all $\{a,b\} \subseteq S$.*

We proceed to prove that AG-groups are, precisely, those strongly regular AG\*\*-groupoids that are left simple.

**Lemma 5 [ 4, Proposition 2.1 and (5), pages 204-205]** *In a completely inverse AG\*\*-groupoid S with $\{a,b\} \subseteq S$ and $\{e, f\} \subseteq S$ ,*
*(1) $a S = S a = a^{-1} S = (aa^{-1}) S$ is the minimal left, right and two-sided ideal containing a,*
*(2) $eS = fS$ implies $e = f$ and*
*(3) e(ab) = (ea)b.*

The next Lemma follows easily from Lemma 5. The proof is omitted.

**Lemma 6.** *In a completely inverse AG\*\*-groupoid S with $\{a,b\} \subseteq S$, the following statements are equivalent:*
*(1) $a S = b S$ ,*
*(2) $aa^{-1} \in b S$ and $bb^{-1} \in a S$*
*(3) $aa^{-1} = bb^{-1}$ .*

**Lemma 7.** *A left simple AG-groupoid S with left identity element $e \in S$ is strongly regular and $E(S) = \{e\}$.*

*Proof* Let $a \in S$. An AG-groupoid with left identity is easily seen to be an AG\*\*-groupoid. So, by Lemma 5, $Sa <_l S$. Therefore $Sa = S$. Hence, there exists $x \in S$ such that

$xa = e = ee = (xa)e = (ea)x = ax$. Therefore, $a = ea = (ax)a$ , with $ax = xa$. By definition, *S* is strongly regular.

Let $f \in E(S)$. By Lemma 5, $Sf <_l S$ and so $Sf = S$. Therefore $xf = e$ for some $x \in S$. Then, $xf = e = ee = (xf)e = (ef)x = fx$ and $f = ef = (xf)f = fx = e$. ∎

**Example 2.** Let $S = \{a,b,c,d,e\}$ be an AG-groupoid with the following Cayley table.

| • | a | b | c | d | e |
|---|---|---|---|---|---|
| a | a | a | a | a | a |
| b | a | b | c | d | e |
| c | a | e | b | c | d |
| d | a | d | e | b | c |
| e | a | c | d | e | b |

By routine calculation it is straightforward to prove that $S$ is a strongly regular AG-groupoid with left identity element $e$. However, $S$ is **not** left simple since $\{a\} = Sa <_l S$.

**Theorem 8.** *A groupoid $S$ is an AG-group if and only if $S$ is a left simple, completely inverse AG\*\*-groupoid.*

*Proof* $(\Rightarrow)$ An AG-group $S$ is an AG-groupoid with a left identity $e$ and so it is an AG\*\*-groupoid. By Lemma 4, for all $a \in S$, $aa^{-1} = a^{-1}a = e$ and so, clearly, $S$ is strongly regular. So, by Theorem 3, $S$ is completely inverse. Suppose that $\emptyset \neq I <_l S$ and let $i \in I$. Then for any $s \in S$, $s = es = (i^{-1}i)s \in IS$. By Lemma 5, $IS = SI \subseteq I$. Thus, $s \in I$ and so $I = S$. So $S$ is left simple.

$(\Leftarrow)$ Let $S$ be a left simple, completely inverse AG\*\*-groupoid. Then by Lemma 5, $S = aS = bS$, for all $\{a,b\} \subseteq S$. It follows from Lemma 6 that $aa^{-1} = bb^{-1}$ $(a,b \in S)$. Then, $(aa^{-1})(sa) = s(aa^{-1} \bullet a) = sa$, for any $s \in S$. Hence, $aa^{-1}$ is the left identity of $S$ and, clearly, $S$ is an AG-group. ∎

**Corollary 9.** *A groupoid $S$ is an AG-group if and only if $S$ is a left simple, AG-groupoid with a left identity element.*

*Proof* $(\Rightarrow)$ It follows from Theorem 8 that an AG-group is left simple. The rest follows from the definition of an AG-group.
$(\Leftarrow)$ This follows from Lemma 7, Theorems 3 and 8 and the fact that an AG-groupoid with a left identity is an AG\*\*-groupoid. ∎

## 5. Inflations of completely inverse AG\*\*-groupoids

**Definition.** An AG-groupoid $S$ is an inflation of its AG-subgroupoid $U$ if $S$ is a disjoint union of the sets $S_u$ $(u \in U)$, $u \in S_u$ $(u \in U)$, and for every $\{a,b\} \subseteq S$ $(a \in S_u, b \in S_v)$, $ab = uv$.

**Theorem 10.** *Let a groupoid $S$ satisfy the medial law. Then $S$ is an inflation of a completely inverse AG\*\*-groupoid if and only if $S^2$ is a completely inverse AG\*\*-groupoid.*

*Proof* $(\Rightarrow)$ Let $S$ be an inflation of its completely inverse, AG\*\*-subgroupoid $T$. Then $T = T^2 \subseteq S^2 \subseteq T^2 = T$ and so $S^2 = T$ is a completely inverse, AG\*\*-groupoid.
$(\Leftarrow)$ Let $\{a,b\} \subseteq S$. Then, since $S$ is medial,

$$ab = (ab)(ab)^{-1} \bullet (ab) = (ab)^{-1}(ab) \bullet (ab) = (ab)^2 (ab)^{-1} = (a^2 b^2)(ab)^{-1}.$$

Therefore,
$$[a^2(a^2)^{-1} \bullet a][b^2(b^2)^{-1} \bullet b] = [a^2(a^2)^{-1} \bullet b^2(b^2)^{-1}](ab) =$$
$$= [a^2(a^2)^{-1} \bullet b^2(b^2)^{-1}][(a^2b^2)(ab)^{-1}]$$

Then, using Lemma 5, part (3) we have
$$[a^2(a^2)^{-1} \bullet a][b^2(b^2)^{-1} \bullet b] = \{[a^2(a^2)^{-1} \bullet b^2(b^2)^{-1}](a^2b^2)\}(ab)^{-1} =$$
$$= \{[a^2(a^2)^{-1} \bullet a^2][b^2(b^2)^{-1} \bullet b^2]\}(ab)^{-1} = (a^2b^2)(ab)^{-1} = ab$$

Define $S_{ab} = \{x \in S : x^2(x^2)^{-1} \bullet x = ab\}$. Then it is straightforward to prove that $S$ is a disjoint union of the sets $S_{ab}$ ($ab \in S^2$), $ab \in S_{ab}$ ($ab \in S^2$) and if $x \in S_{ab}$ and $y \in S_{cd}$ ($ab, cd \in S^2$) then $xy = ab \bullet cd$. Hence, $S$ is an inflation of $S^2$. ∎

## 6. Completely inverse AG**-groupoids and commutative inverse semigroups

A semigroup is commutative and inverse if and only if it is a semilattice of abelian groups. It is well known that a semigroup $S$ is an inflation of a semilattice of groups if and only if $S^2$ is a semilattice of groups. Theorem 10 then hints that completely inverse $AG^{**}$-groupoids are closely related to semilattices of abelian groups and we proceed to elaborate this relationship.

**Definition.** If $(S, \bullet)$ is a completely inverse $AG^{**}$-groupoid then we define a product $[\bullet]$ on $S$ as $a[\bullet]b = (b \bullet bb^{-1})a$, for $\{a,b\} \subseteq S$.

**Proposition 11.** *If $(S, \bullet)$ is a completely inverse $AG^{**}$-groupoid then the identity mapping id on $E(S, \bullet)$ is an isomorphism between $E(S, \bullet)$ and $E(S, [\bullet])$.*

*Proof* If $e \in E(S, \bullet)$ then $e[\bullet]e = e \bullet e = e \in E(S, [\bullet])$. Conversely, if $x \in E(S, [\bullet])$ then
$x[\bullet]x = x = (x \bullet xx^{-1})x$. So $xx^{-1} = [(x \bullet xx^{-1})x]x^{-1} = (xx^{-1})(x \bullet xx^{-1}) = x \bullet (xx^{-1})^2 =$
$= x \bullet xx^{-1} = (xx^{-1})^2 = (x \bullet xx^{-1})(xx^{-1}) = (xx^{-1})^2 \bullet x = xx^{-1} \bullet x = x = x \bullet x \in E(S, \bullet)$.
So we have proved that $E(S, \bullet) = E(S, [\bullet])$. Then for $e, f \in E(S, \bullet)$,
$id(e \bullet f) = e \bullet f = f \bullet e = (f \bullet ff^{-1}) \bullet e = e[\bullet]f = (id\,e)[\bullet](id\,f)$. Hence,
$id: E(S, \bullet) \to E(S, [\bullet])$ is an isomorphism. ∎

**Theorem 12.** *If $(S, \bullet)$ is a completely inverse $AG^{**}$-groupoid then $[\bullet]$ is commutative and associative and $(S, [\bullet])$ is a semilattice of abelian groups.*

*Proof* Let $\{a,b,c\} \subseteq S$. Then,
$a[\bullet]b = (b \bullet bb^{-1})a = (b \bullet bb^{-1})(aa^{-1} \bullet a) = (aa^{-1})[(b \bullet bb^{-1})a] = (aa^{-1})[(a \bullet bb^{-1})b] =$
$= (a \bullet bb^{-1})(aa^{-1} \bullet b) = (a \bullet aa^{-1})(bb^{-1} \bullet b) = (a \bullet aa^{-1})b = b[\bullet]a$.
Hence, $(S, [\bullet])$ is a commutative groupoid. Also, note that $(a[\bullet]b)(a[\bullet]b)^{-1} =$
$= [(b \bullet bb^{-1})a][(b \bullet bb^{-1})a]^{-1} = [(b \bullet bb^{-1})a][(b^{-1} \bullet bb^{-1})a^{-1}] = (bb^{-1})(aa^{-1})$.
Then, using Proposition 11, the idempotents of $(S, \bullet)$ form a semilattice. Therefore, using Lemma 5, part (3), we have
$a[\bullet](b[\bullet]c) = \{(b[\bullet]c)[(b[\bullet]c)(b[\bullet]c)^{-1}]\}a = \{[(c \bullet cc^{-1})b][(cc^{-1})(bb^{-1})]\}a =$
$= \{[(c \bullet cc^{-1})](cc^{-1})](b \bullet bb^{-1})\}a = [c(b \bullet bb^{-1})]a = [(cc^{-1} \bullet c)(b \bullet bb^{-1})]a =$
$= \{(cc^{-1})[c(b \bullet bb^{-1})]\}a = (cc^{-1})\{[c(b \bullet bb^{-1})]a\} = [c(b \bullet bb^{-1})](cc^{-1} \bullet a) =$
$= (c \bullet cc^{-1})[(b \bullet bb^{-1})a] = (a[\bullet]b)[\bullet]c$.

Therefore $(S,\bullet)$ is a commutative semigroup. By Proposition 11, $E(S,\bullet)$ is a semilattice. It is easy to calculate that $(a^{-1}\bullet aa^{-1})[\bullet]a[\bullet]a = a$. It follows that $a \in a^2[\bullet]S = S[\bullet]a^2$ and so by [2, Theorem 4.3 and Theorem 4.11] $(S,[\bullet])$ is a semilattice of abelian groups. ∎

**Corollary 13.** *If $(S,\bullet)$ is a completely inverse AG\*\*-groupoid and $a \in S$ then $a \bullet a^{-1}{}^{-1} = a^{-1} \bullet aa^{-1}$ is the inverse of $a$ in $(S,[\bullet])$ and $a \bullet aa^{-1} = a^{-1} \bullet aa$ is the inverse of $a^{-1}$ in $(S,[\bullet])$. Also, in $(S,[\bullet])$, $aa^{-1}$ is the identity of the group to which $a$ belongs.*

**Notation.** If $(S,\bullet)$ is a completely inverse $AG^{**}$-groupoid and $a \in S$ then $a_\bullet^{-1}$ or $(a)_\bullet^{-1}$ denotes $a^{-1}$, the inverse of $a$ in $(S,\bullet)$. Therefore, $a_{[\bullet]}^{-1} = (a)_{[\bullet]}^{-1} = a \bullet a^{-1}{}^{-1} = a^{-1} \bullet aa^{-1}$ and $(a^{-1})_{[\bullet]}^{-1} = a \bullet aa^{-1} = a^{-1} \bullet aa$.

**Proposition 14.** *If $(S,\bullet)$ is a semilattice of abelian groups then $(S,\bullet) \cong (S,[\bullet])$.*

*Proof* Let $\{a,b\} \subseteq S$. Then $a[\bullet]b = (b \bullet bb^{-1})a = (bb^{-1} \bullet b)a = b \bullet a = a \bullet b$. ∎

**Theorem 15.** *If $(S,\bullet)$ and $(T,\circ)$ are completely inverse $AG^{**}$-groupoids then $(S,[\bullet]) \cong (T,[\circ])$ if and only if there is a bijection $\boldsymbol{B}:(S,\bullet) \to (T,\circ)$ such that for every $\{a,b\} \subseteq S$, $\boldsymbol{B}(a \bullet b) = [\boldsymbol{B}(a)_\bullet^{-1}]_\circ^{-1} \circ \boldsymbol{B}b$.*

*Proof* $(\Rightarrow)$ Suppose that $\boldsymbol{B}:(S,[\bullet]) \cong (T,[\circ])$. Note that, using the medial and paramedial laws,
$$a \bullet b = (aa^{-1} \bullet a)(bb^{-1} \bullet b) = (b \bullet a)(bb^{-1} \bullet aa^{-1}) = (b \bullet bb^{-1})(a \bullet aa^{-1}).$$
Therefore,
(1) $a \bullet b = (a \bullet aa^{-1})[\bullet]b$ and similarly
(2) $\boldsymbol{B}a \circ \boldsymbol{B}b = [\boldsymbol{B}a \circ \boldsymbol{B}a\,(\boldsymbol{B}a)^{-1}][\circ]\boldsymbol{B}b = [\boldsymbol{B}b \circ \boldsymbol{B}b\,(\boldsymbol{B}b)^{-1}][\boldsymbol{B}a \circ \boldsymbol{B}a\,(\boldsymbol{B}a)^{-1}]$.

Since $\boldsymbol{B}:(S,[\bullet]) \cong (T,[\circ])$, and using (1), we have
(3) $\boldsymbol{B}(a \bullet b) = \boldsymbol{B}[(a \bullet aa^{-1})[\bullet]b] = \boldsymbol{B}(a \bullet aa^{-1})[\circ]\boldsymbol{B}b$

But from Corollary 13, we have
(4) $\boldsymbol{B}(a \bullet aa^{-1}) = \boldsymbol{B}[(a)_\bullet^{-1}]_{[\bullet]}^{-1} = [\boldsymbol{B}(a)_\bullet^{-1}]_{[\circ]}^{-1} = [\boldsymbol{B}(a)_\bullet^{-1}]_\circ^{-1} \circ \boldsymbol{B}[(a)_\bullet^{-1}]\{\boldsymbol{B}[(a)_\bullet^{-1}]\}_\circ^{-1}$

Therefore, using (3), (4) and paramediality we have
$\boldsymbol{B}(a \bullet b) = \boldsymbol{B}(a \bullet aa^{-1})[\circ]\boldsymbol{B}b = [\boldsymbol{B}b \circ \boldsymbol{B}b(\boldsymbol{B}b)^{-1}] \circ \boldsymbol{B}(a \bullet aa^{-1}) =$
$= [\boldsymbol{B}b \circ \boldsymbol{B}b(\boldsymbol{B}b)^{-1}] \circ \{[\boldsymbol{B}(a)_\bullet^{-1}]_\circ^{-1} \circ \boldsymbol{B}(a)_\bullet^{-1}[\boldsymbol{B}(a)_\bullet^{-1}]_\circ^{-1}\} = [\boldsymbol{B}(a)_\bullet^{-1}]_\circ^{-1} \circ \boldsymbol{B}b$

$(\Leftarrow)$ Suppose that $\boldsymbol{B}:(S,\bullet) \to (T,\circ)$ is a bijection such that for $\{a,b\} \subseteq S$, we have $\boldsymbol{B}(a \bullet b) = [\boldsymbol{B}(a)_\bullet^{-1}]_\circ^{-1} \circ \boldsymbol{B}b$. Then $\boldsymbol{B}(a[\bullet]b) = \boldsymbol{B}[(b \bullet bb^{-1})a] = [\boldsymbol{B}(b^{-1} \bullet bb^{-1})]_\circ^{-1} \circ \boldsymbol{B}a =$
$= [(\boldsymbol{B}b)_\circ^{-1} \circ \boldsymbol{B}b(\boldsymbol{B}b)_\circ^{-1}]_\circ^{-1} \circ \boldsymbol{B}a = [(\boldsymbol{B}b) \circ \boldsymbol{B}b(\boldsymbol{B}b)_\circ^{-1}]_\circ^{-1} \circ \boldsymbol{B}a = \boldsymbol{B}a[\circ]\boldsymbol{B}b$
and therefore $\boldsymbol{B}:(S,[\bullet]) \cong (T,[\circ])$ ∎

**Corollary 16.** *If $(S,\bullet)$ is a completely inverse AG\*\*-groupoid then a bijection $\boldsymbol{B}: S \to S$ is an automorphism on $(S,[\bullet])$ if and only if for any $a,b \in S$, $\boldsymbol{B}(a \bullet b) = [\boldsymbol{B}(a_\bullet^{-1})]_\bullet^{-1} \bullet \boldsymbol{B}b$.*

*Proof* This follows from Proposition 14 and the proof of Theorem 15, with $(T,\circ) = (S,[\bullet])$. ∎

We now define a relation $\approx$ on the class of completely inverse $AG^{**}$-groupoids as follows: $(S,\bullet) \approx (T,\circ)$ if and only if $(S,[\bullet]) \cong (T,[\circ])$. It is then obvious that:

**Proposition 17.** *The relation $\approx$ is an equivalence relation.*

It is also clear from Proposition 14 that:

**Proposition 18.** *If $(S,\bullet)$ is a completely inverse AG\*\*-groupoid then, up to isomorphism, the $\approx$ equivalence class containing $(S,\bullet)$ contains exactly one semilattice of abelian groups; in particular, $(S,[\bullet]) \approx (S,\bullet)$.*

We define a mapping $\Omega$ from the class of semilattices of abelian groups to the collection of $\approx$ equivalence classes as follows: $\Omega(S,\bullet)$ is defined as the $\approx$ equivalence class containing $(S,\bullet)$. If we consider two semilattices of abelian groups to be equal if they are isomorphic, then it follows readily from Proposition 17 that:

**Proposition 19.** *The mapping $\Omega$ is a bijection and if $(S,\bullet),(T,\circ)$ are two semilattices of abelian groups then $(S,\bullet) \cong (T,\circ)$ if and only if $\Omega(S,\bullet) = \Omega(T,\circ)$.*

## 7. The structure of completely inverse AG\*\*-groupoids

It seems from Theorems 15 and 19 in Section 6 that the structure of a completely inverse AG\*\*-groupoid $(S,\bullet)$ is closely related to the structure of the commutative, inverse semigroup in its $\approx$ class; namely, $(S,[\bullet])$. However, there is no obvious way to use those theorems to construct the product $\bullet$ from that of $[\bullet]$.

In this section, building on the work of Sharma [7] and Shah, Shpectorov and Ali [6], we prove that all completely inverse AG\*\*-groupoids can be constructed from their $\approx$ representative $(S,[\bullet])$. To do this we need to know the product $[\bullet]$ and all involutive, E-fixed automorphisms on $(S,[\bullet])$.

**Theorem 20.** *Let $(S,\bullet)$ be a semilattice of abelian groups and $A \in AUT_e^2(S)$. Define a product $\circ$ on $S$ as follows: $a \circ b = Aa \bullet b$ $(a,b \in S)$. Then $(S,\circ)$ is a completely inverse AG\*\*-groupoid and $(S,\bullet) \cong (S,[\circ])$.*

*Proof* Let $\{a,b,c\} \subseteq S$. Then $(a \circ b) \circ c = (Aa \bullet b) \circ c = [A(Aa \bullet b)] \bullet c = A^2 a \bullet Ab \bullet c =$
$= a \bullet Ab \bullet c = c \bullet Ab \bullet a = (c \circ b) \circ a$. Also, $a \circ (b \circ c) = Aa \bullet Ab \bullet c = Ab \bullet Aa \bullet c = b \circ (a \circ c)$.
So $(S,\circ)$ is an AG\*\*-groupoid.

Also, $(a \circ Aa^{-1}) \circ a = A(Aa \bullet Aa^{-1}) \bullet a = (A^2 a \bullet A^2 a^{-1}) \bullet a = (a \bullet a^{-1}) \bullet a = a$ and
$(Aa^{-1} \circ a) \circ Aa^{-1} = (A^2 a^{-1} \bullet a) \circ Aa^{-1} = A(a^{-1} \bullet a) \bullet Aa^{-1} = A(a^{-1}a \bullet a^{-1}) = Aa^{-1}$
Hence, $Aa^{-1}$ is an inverse of $a$.
As noted in the preliminaries, since $(S,\circ)$ is an AG\*\*-groupoid, $Aa^{-1}$ is the unique inverse of $a$.
Then, $Aa^{-1} \circ a = A^2 a^{-1} \bullet a = a^{-1}a = aa^{-1} = A(aa^{-1}) = Aa \bullet Aa^{-1} = a \circ Aa^{-1}$ and so $(S,\circ)$ is a completely inverse AG\*\*-groupoid. Also, $a[\circ]b = [b \circ (b \circ b_\circ^{-1})] \circ a = [Ab \bullet (Ab \bullet Ab^{-1})] \circ a =$
$= Ab \circ a = b \bullet a = a \bullet b$ and so $(S,\bullet) \cong (S,[\circ])$. ∎

**Theorem 21.** *If $(S,\bullet)$ is a completely inverse AG\*\*-groupoid, define a mapping $A: S \to S$ by $Aa = a \bullet aa^{-1}$. Then $A \in AUT_e^2(S,[\bullet])$ and $a \bullet b = Aa[\bullet]b$.*

*Proof* Let $\{a,b\} \subseteq S$. Then $A(a[\bullet]b) = A[(b \bullet bb^{-1})a] = [(b \bullet bb^{-1})a] \bullet (bb^{-1} \bullet aa^{-1}) =$
$= [(b \bullet bb^{-1})(bb^{-1})](a \bullet aa^{-1}) = (a \bullet aa^{-1})[\bullet](b \bullet bb^{-1}) = Aa[\bullet]Ab$.
So, $A$ is a homomorphism on $(S,[\bullet])$.

If $a \bullet aa^{-1} = b \bullet bb^{-1}$ then $b = (b \bullet bb^{-1})(bb^{-1}) = (a \bullet aa^{-1})(bb^{-1}) = (bb^{-1} \bullet aa^{-1}) \bullet a =$
$= (bb^{-1})(aa^{-1} \bullet a) = (bb^{-1})a$. Dually, $a = (aa^{-1})b$. Then, $a = (aa^{-1})b = (aa^{-1})(bb^{-1} \bullet a) =$
$= (bb^{-1})(aa^{-1} \bullet a) = (bb^{-1})a = b$ and so $A$ is one-to-one.

Since $(a \bullet aa^{-1}) \bullet (a \bullet aa^{-1})(a \bullet aa^{-1})^{-1} = (a \bullet aa^{-1}) \bullet aa^{-1} = a$, $A$ is an involutive automorphism and, clearly, $A \in AUT_e^2(S,[\bullet])$. Finally, using the paramedial law,
$Aa[\bullet]b = (b \bullet bb^{-1})(a \bullet aa^{-1}) = a \bullet b.$ ∎

**Corollary 22** *A groupoid* $(S, \bullet)$ *is a completely inverse AG\*\*-groupoid if and only if there exists a product* $\circ$ *on S, such that* $(S, \circ)$ *is a semilattice of groups, and an* $A \in AUT_e^2(S, \circ)$ *such that* $a \bullet b = Aa \circ b$ $(a,b \in S).$

[1] 10 Albert Mansions, London N8 9RE, United Kingdom, e-mail: bobmonzo@talktalk.net